\documentclass[12pt]{amsart}
\usepackage{amsthm,amsmath,amssymb,amscd,graphics,enumerate}
\usepackage[all]{xy}

\newcommand{\thisdate}{\today}
{
   
   \newtheorem{proposition}[subsubsection]{Proposition}     
   \newtheorem{lemma}[subsubsection]{Lemma}

   \newtheorem{corollary}[subsubsection]{Corollary}

}
{\theoremstyle{definition}

   \newtheorem{definition}[subsubsection]{Definition}
   
}
\newcommand{\RR}{{\mathbb{R}}}

\newcommand{\HH}{{\mathbb{H}}}
\newcommand{\CC}{{\mathbb{C}}}

\newcommand{\PP}{{\mathbb{P}}}
\newcommand{\ZZ}{{\mathbb{Z}}}

\newcommand{\GG}{{\mathbb{G}}}

\newcommand{\bbA}{{\mathbb{A}}}

\newcommand{\bQ}{{\mathbf{Q}}}

\newcommand{\m}{{\mathfrak{m}}}

\newcommand{\cI}{{\mathcal I}}

\newcommand{\cM}{{\mathcal M}}

\newcommand{\cO}{{\mathcal O}}
\newcommand{\cP}{{\mathcal P}}

\newcommand{\cV}{{\mathcal V}}
\newcommand{\cX}{{\mathcal X}}
\newcommand{\cY}{{\mathcal Y}}
\newcommand{\cZ}{{\mathcal Z}}

\newcommand{\Spec}{\operatorname{Spec}}
\newcommand{\Rank}{\operatorname{Rank}}

\newcommand{\cProj}{{\cP}roj}

\newcommand{\Pic}{{\operatorname{\mathbf{Pic}}}}

\newcommand{\lrar}{\longrightarrow}
\newcommand{\dar}{\downarrow}

\newcommand{\ocZ}{\overline{{\mathcal Z}}}
\newcommand{\double}{\genfrac..{0pt}1
{\raise -1pt\hbox{$\scriptstyle\longrightarrow$}}{\raise 3pt\hbox
{$\scriptstyle\longrightarrow$}}} 

\renewcommand{\setminus}{\smallsetminus}
\newcommand{\setmin}{\smallsetminus}

\def\tototi{\mathbin{\mathop{\otimes}\limits^{\raise-1pt\hbox
{$\scriptscriptstyle {\rm L}$}}}}

\def\indlim{\mathop{\vrule width0pt height7pt depth
4pt\smash{\lim\limits_{\raise 1pt\hbox to 14.5pt
{\rightarrowfill}}}}}
\def\projlim{\mathop{\vrule width0pt height7pt depth
4pt\smash{\lim\limits_{\raise 1pt\hbox to 14.5pt
{\leftarrowfill}}}}}

\begin{document}
\hfill \date{\thisdate}

\title{Computations with moduli of perverse point sheaves}  
\author[D. Abramovich]{Dan Abramovich}
\thanks{Research of D.A. partially supported by NSF grant DMS-0070970}  
\address{Department of Mathematics\\ Boston University\\ 111 Cummington
         Street\\ Boston, MA 02215\\ U.S.A.} 
\curraddr{Department of Mathematics,
Box 1917, 
Brown University, Providence, RI, 02912}
\email{abrmovic@math.bu.edu,abrmovic@math.brown.edu}
\author{Jiun C.\ Chen}
\address{Department of Mathematics\\ Harvard  University\\ 1 Oxford
         Street\\ Cambridge, MA 02138\\ U.S.A.} 
\email{jcchen@math.harvard.edu}
\maketitle

\section{Introduction}
We work with schemes of finite type over $\CC$.

\subsection{Setup and Bridgeland's results}
Let $f: X \to Y$ be a projective birational contraction satisfying
the following two assumptions: 
\begin{enumerate}
\item[\bf (B1) ] $\RR f_* \cO_{X} = \cO_Y$, and
\item[\bf (B2) ] $\dim f^{-1} \{z\} \leq 1$ for every $z\in Y$.
\end{enumerate}
 In \cite{B}, T.\ Bridgeland constructs a moduli space $\cM(X/Y)$ of
 {\em perverse point sheaves} on $X/Y$. This moduli space admits a
 natural morphism $\cM(X/Y)\to Y$, and there is an
 irreducible component $W(X/Y)\subset \cM(X/Y)$ - the main
 component of $\cM(X/Y)$ -  such that 
 $W(X/Y)\to Y$ is birational. 

In the main result of \cite{B}, one further assumes
 that $X$ is a smooth threefold and $X\to Y$ is a
flopping contraction. In this situation it is shown that
 $W(X/Y)\to Y$ is isomorphic to the 
flop $X^+ \to Y$, and further that the correspondence given by the
universal perverse point sheaf gives an equivalence of derived 
categories $D_{coh}( X^+ ) \to D_{coh}( X )$. This result was generalized by
 the second author in \cite{Chen} to the case where $X$ is a terminal
 Gorenstein threefold. Further results on equivalence of derived categories
 were proven by Bondal and Orlov \cite{Bondal-Orlov}, Kawamata,
 e.g. \cite{Kawamata}  and others.

In this note we collect some results orthogonal to the derived category
discussion about Bridgeland's moduli space of perverse point sheaves. Some of
these are used in the paper \cite{AC}, to which this note is a companion.

\subsection{The main result} Let $f:X \to Y$ be a projective
 birational contraction, and let 
 $M$ be an $f$-ample line bundle on $X$. If the sheaf of
 $\cO_Y$-algebras 
$$A^+ := \mathop\bigoplus\limits_{i\leq 0}\ f_* M^i$$
 is finitely
 generated, we denote 
$$X^+ = \cProj_{_Y} A^+.$$
 As part of the minimal
 model program, it is conjectured that $X^+ \to Y$ exists, in other
 words, 
 $A^+$ is finitely generated, whenever $f$ is a flipping or flopping
 contraction (or, more generally, a log-flipping contraction with
 respect to a suitable divisor); in such situation $X^+ \to Y$ is
 called the flip, respectively  flop, of $X \to Y$.  

Now assume $f:X \to Y$ satisfies conditions (B1)-(B2) above. It is
tempting to ask if $W(X/Y)$ is always isomorphic to $X^+$, 
or, at least, for criteria for comparing $W(X/Y)$ with $X^+$. It
turns out that the {\em normalization} $W^n$ of $W(X/Y)$ is
more approachable. 

The first result concentrates on the case when $X\to Y$ and $W \to
Y$ do not contract divisors:

\begin{proposition}\label{Prop:proj}
Let $f:X\to Y$ be a birational morphism satisfying (B1) and (B2)
above.  Let $M$ be an $f$-ample invertible sheaf. Denote by $W$ the
main component of the moduli space of perverse 
point sheaves on $X/Y$, and let $W^n \to W$ be the normalization.

Assume further that
  $X$ and $Y$ are {\em normal} algebraic
verieties and the exceptional loci of $X\to Y$ and $W \to Y$ have codimension
$>1$. Then  $$W^n = \cProj_Y\oplus_{i\leq 0} M^{i}.$$
\end{proposition}


We denote by $U$ the open set where $X \to Y$ is an isomorphism.
The following {\em flatness criterion} has been applied
in \cite{Chen}:

\begin{proposition}\label{Prop:flatness}
Let $f:X \to Y$ be a projective birational morphism of normal
varieties, satisfying (B1)-(B2) above. 

Let  $W \to Y$ be the main component of the moduli space of perverse
point sheaves on $X/Y$. Then
\begin{enumerate}
\item the ideal sheaf $\cI_{X\times_Y W}$ of the fibered product
$X\times_Y W$ inside $X \times W$ is flat over $W$,
\item the universal perverse point sheaf on $X \times W$ is the
structure sheaf of the fibered product $X\times_Y W$, and
\item if $V$ is a scheme and $h:V \to Y$ is a morphism such that every
associated component of $V$ intersects $h^{-1} U$, then $h$
factors through $W \to Y$ if and only if the ideal  $\cI_{X\times_Y
V}$ of the fibered product $X\times_Y V$ in $X \times V$ is flat over $V$.
\end{enumerate}
\end{proposition}

Combining the two propositions we obtain:

\begin{corollary}\label{Cor:flp} let $X \to Y$ be a small birational
contraction. Assume $X^+$ exists and  the ideal sheaf of $X 
\times_Y X^+$ is flat over  $X^+$. Then $W^n = X^+$.
\end{corollary}

\subsection{Computations} 
To illustrate the possible applications of our propositions we 
use our results to test some basic cases of birational contractions.

The first case of interest is that of a surface contraction $X \to
Y$. This case does not fall under Proposition \ref{Prop:proj}, but we
can still use the flatness criterion\ref{Prop:flatness} to study the 
situation.  

We have the following basic cases:

\begin{proposition}\label{Prop:surfaces}
\begin{enumerate}
\item Suppose $X\to Y$ is a proper birational map between smooth surfaces. Then
$W = Y$.   
\item Suppose $X\to Y$ is a contraction of an irreducible (-2)-curve
on a smooth surface
$X$. Then $W^n$ is isomorphic to $X$.
\end{enumerate}
\end{proposition}

We remark that in case (2) above, and indeed for an arbitrary crepant
contraction of surfaces with canonical singularities, the results of
\cite{Chen} imply that $W$ is isomorphic to $X$.

The first case of a threefold flip is described in \cite{Kollar-Mori},
page 39: it is the Francia flip $X_n \to Y_n \leftarrow X_n^+$ of
index $n$. 

We have:

\begin{proposition}\label{Prop:Francia}
\begin{enumerate}
\item $W^n(X_n/Y_n) = X_n^+$. 
\item The morphism $X_2\to Y_2$ does not factor through the moduli space
$W(X_2^+/Y_2)$.  
\end{enumerate}
\end{proposition}

Part (2) is a natural example of a log-flip (i.e. the inverse of the Francia
flip) which is not directly related to 
Bridgeland's construction. It should not be difficult to show that the same
holds when $2$ is replaced by $n>1$.

In higher dimensions we have:

\begin{proposition}\label{Prop:segre}
Let $Y$ be the cone over the Segre embedding of $\PP^1 \times
\PP^k$. Let $X\to Y$ be the small resolution having $\PP^1$ as the
exceptional locus. Then $W^n(X/Y)  = X^+$, the other small resolution.
\end{proposition}

\subsection{Related work} 
\begin{enumerate}
\item In \cite{Chen}, the second author treats the
case of a terminal treefold flopping contraction $X \to
Y$ with hypersurface singularities. In that case the full analogue of
Bridgeland's theorem is shown.
\item In \cite{ABC}, the Francia flip of index $n$ is treated in
detail. It is shown that by endowing $X^-$ with a natural smooth stack
structure $\cX^-$, and
using perverse point sheaves of perversity $(-1,0)$ instead of the
definition of \cite{B}, one gets an isomorphism $W^{(-1,0)}(\cX^-x/Y) =
X^+$ and a fully faithful functor on derived catogories $D_{coh}(X^+)
\to D_{coh}(\cX)$.
\item In \cite{Kawamata} the general case of a toric flip is
considered without moduli spaces. It is shown that, again using 
natural orbifold (or stack) structures, that, in a wide variety of
cases one obtains a fully faithful functor of derived categories.
\end{enumerate}

\subsection{Acknowledgements} We thank Tom Bridgeland, David Eisenbud
and Johan de Jong for helpful discussions.

\section{Proof of Proposition \ref{Prop:proj}}

\subsection{Perverse point sheaves and perverse ideal sheaves}\label{Sec:PPS}
 In \cite{B}, the moduli space $\cM(X/Y)$ is defined via an exotic
 $t$-structure on $D_{coh}(X)$, whose heart is
 the category of perverse coherent sheaves $Per(X/Y)$.
 However,  for  what we say below   we only need the following:
\begin{definition}
 A {\em perverse point sheaf} (PPS) for $\cX/Y$ is a complex $E = [F \to \cO_\cX]$,
 with $\cO$ positioned in degree 0, where
\begin{enumerate}
\item the sheaf $f_*(F)$ is an ideal sheaf of a point in $Y$, with an
  embedding $f_*(F) \to \cO_Y$ induced by pushforward from $E = [F \to
  \cO]$, and
\item the natural map $f^*f_*(F) \to F$ is surjective.
\end{enumerate}
The sheaf $F$ is called {\em the perverse point-ideal sheaf} (PPIS)
associated to the PPS $E$.
\end{definition}

We note that, by the naturality of $f^*$ and the fact that $f^*\cO_Y =
\cO_\cX$ we have that the composition of $f^*f_*(F) \to F$ with the
given $F \to \cO_\cX$ coincides with the map  $f^*f_*(F) \to \cO_\cX$
obtained by pulling back $f_*(F) \to \cO_Y$.

\subsection{Brief account of Bridgeland's construction}
In \cite{B}, section 5, the moduli space of perverse point sheaves is constructed
as the moduli space of the corresponding perverse {\em point-ideal} sheaves
using geometric 
invariant theory. Let $U\subset X$ be the open set where $f: X \to Y$
is an  isomorphism, and pick a point $x \in U$. The maximal ideal
sheaf $F_0 = I_x$ is a 
perverse point-ideal sheaf associated to the perverse point sheaf
$\cO_x$. Pick an ample line bundle $L_Y$ on $Y$, and write $L= f^*
L_Y$. One considers the vector space  
$$ V =  Hom (L_Y^{-1}, f_* F_0) = H^0(Y, L_Y \otimes f_* F_0).$$

It is shown that every perverse point-ideal sheaf $F$ is a quotient of
$V \otimes L$.  Such a quotient is characterised by the surjective linear map
$$H^0(X, V \otimes L \otimes M^n) \to H^0(X, F \otimes M^n),$$ for a
sufficiently large integer $n$. Such a map corresponds to a point is a
Grassmannian $\GG$. It is shown that, for large $n$, the moduli space
is the GIT quotient of 
the set of semistable points $\cZ$ on a closed subset $\ocZ \subset
\GG$, linearised 
through the tautological line bundle $\cO_{\PP^N}(1)$ of the  Pl\"ucker
embedding $\GG \subset \PP^N$. All these semistable points  are strictly
stable. 

\subsection{Comparison of line bundles}

Denote by $(\PP^N)^s$ the locus of semistagble points, and similarly
$\GG^s\subset \GG$. 
We use the following diagram:

$$\begin{array}{cccccc}
\ocZ & \subset & \GG   & \subset &\PP^N  \\
 \cup&         &\cup   &         &\cup      \\
\cZ  & \subset & \GG^s & \subset &(\PP^N)^s \\
\dar &         &       &         & \dar \\
 \cM(X/Y) &    &\stackrel{i}{\hookrightarrow}&         & \bQ \\
\end{array}$$
where $\bQ$ is the projective spectrum of the ring of invariant
sections on $\PP^N$. By   
construction we have that the pullback of $\cO_\bQ(1)$ to $(\PP^N)^s $
is the same as the restriction of $\cO_{\PP^N}(1)$.

We  keep the notation $U\subset X$ for the open set where $f: X \to Y$ is
an isomorphism. There is a natural embedding $j :U \to
\cM(X/Y)$.

\begin{lemma} 
 Consider the composite morphism 
$$h \ = \ i\circ j  : \  U \to \bQ.$$

Then there is an integer $k$ and an isomorphism 
$$ h^* \cO_\bQ(1) \simeq j^* ( L^k \otimes M^{-n}).$$ 
\end{lemma}

{\em Proof.} The universal perverse ideal sheaf on $U$ is the ideal
sheaf of the graph $\Gamma_j$ of $j: U \to X$ in $U \times
X$. Denote the two projections by $\pi_U: U \times
X \to U$ and $\pi_X : U \times
X \to X$.  Since $L_Y$ is sufficiently ample and $f$ is isomorphic
over $U$, we have an
exact sequence
$$ 0 \to \ \ \pi_{U\,*} ( I_{\Gamma_g} \otimes \pi_X^* L)\ \  \to\ \ 
\cO_U\otimes H^0(X, L)\ \  \to\ \  j^* L\ \  \to 0.$$
Denote the sheaf on the left hand side of the sequence
$$\cV \ :=\ \ \pi_{U\,*} ( I_{\Gamma_g} \otimes \pi_X^* L)$$ 
and its rank by $\m = \Rank \cV$. It follows that $c_1(\cV) = j^*L^{-1}$.

Let $\HH$ be the total space of 
$Hom (\cO_U^\m,\cV)$ and $g: \HH \to U$ the natural projection. There
is a tautological determinant section $s\in H^0(\HH, g^* L^{-1})$, whose
non-zero locus is the principal frame bundle of $\cV$:
$$\cP\ :=\ \ \HH\setminus div(s)\ \ \ \subset \ \ \ \HH,$$ on which $\cV$
is trivialized. Note  the
isomorphism of  
invertible sheaves
 $\cO_\HH(div(s)) \simeq g^*L^{-1}$.

 The morphism $\cP \to \bQ$ lifts to a
morphism $\cP \to \GG\to \PP^N$ through the quotient
$$ g^*(\cV \otimes \pi_{U\, *} \pi_X^* (L \otimes M^n))  \ \lrar \ 
   g^*\pi_{U \, *} ( I_{\Gamma_j}\otimes \pi_X^* M^n )). $$
The pullback of $\cO_{\PP^N}(1)$ is the determinant of $ g^* \pi_{U \,
*} ( I_{\Gamma_j}\otimes \pi_X^* M^n ).$ Using the exact sequence 
$$ 0 \to\ \ \pi_{U\,*} ( I_{\Gamma_j} \otimes \pi_X^* M^n)\ \ \to\ \
\cO_U\otimes H^0(X, M^n)\ \ \to\ \ j^* M\ \ \to 0$$
we see that  $$g^* h^* \cO_{\PP^N}(1)\ \ \simeq\ \ g^* j^*M^{-n}.$$
 
Since $\HH\to U$  is an affine bundle we have an isomorphism
$g^*:\Pic(U)  \stackrel{\sim}{\to} \Pic(\HH)$, and the excision
sequence gives  
$$\ Pic(\cP) \ \ \ = \ \ \  \Pic(\HH)\ /\ \ZZ\cdot [\cO_\HH(div(s))] \
\ =\ \
  \Pic(U)\ /\ \ZZ\cdot [j^*L].$$ 
The Lemma follows.

\qed

\subsection{Conclusion of proof} 

Since $X \setmin U$ has codimension $\geq 2$ and $X$ is normal, every
section of $ (M^{-n}\otimes L^k )^i|_U$ extends to a section on $X$,
and similarly for $W^n$. The proposition follows. \qed

\section{The flatness criterion}

\subsection{Torsion-freeness of flat sheaves}
We say that an open subscheme $U\subset W$ is schematically dense if
 it contains all  associated primes of W.
 We use the following well-known lemma:

\begin{lemma}
Let $\pi:Z \to W$ be a morphism of schemes,
and let $Q$ be a coherent sheaf on 
$Z$, flat over $W$.   Assume
there is a  schematically dense  open subscheme $U \subset W$ 
such that $Q_{\pi^{-1}U}$ is torsion-free. Then $Q$ is torsion-free. 
\end{lemma} 

{\em Proof.} We argue by contradiction. The question is local, so we
may replace $Z$ and $W$ by 
affine schemes. If $s$ is a torsion section of $Q$, then the image
$T\subset W$
of
its support  is disjoint from the dense open $U$. Let $f$ be a 
function on $W$ which vanishes on $T$. Since $U$ is schematically
dense, $f$ is not a zero divisor in $\cO_W$, thus the sequence  
$$ 0 \to \cO_W \stackrel{f}{\to} \cO_W \to \cO_{V(f)} \to 0$$
is exact. Since $Q$ is flat we have that the sequence 
$$ 0 \to Q \stackrel{f}{\to} Q \to Q_{V(f)} \to 0$$
is exact. But $s$ is in the kernel of the homomorphism $Q
\stackrel{f}{\to} Q$, which is a contradiction.\qed

\subsection{Proof of the proposition}

Let $[F\to \cO_{W\times X}]$ be the universal perverse point
sheaf. The sheaf $F$ is flat over $W$ by definition. Over the open set
$U$, this coincides with the natural 
embedding $[I_{\Gamma_j} \to \cO_{U\times X}]$, and the sheaf
$I_{\Gamma_j}$ is clearly torsion free. By the lemma $F$ is torsion
free. Also, since $[I_{\Gamma_j} \to \cO_{U\times X}]$ is injective,
the kernel of $F\to \cO_{W\times X}$ is torsion, therefore it
vanishes. Thus $F\to \cO_{W\times X}$ is injective and $F$ is an ideal
sheaf. The pushforward $(id \times f)_* F \hookrightarrow \cO_{W
\times Y}$ is the ideal sheaf of the graph of $W \to X$, thus the
image of the pullback
 $$(id \times f)^*(id \times f)_* F \to  \cO_{W\times X}$$
is the ideal sheaf of $W \times_Y X$. But, as we have seen in \ref{Sec:PPS}, this
sheaf homomorphism factors through  $F\to \cO_{W\times X}$. Since the
latter is injective we have that $F$ coincides with the the ideal
sheaf of $W \times_Y X$, which proves (1). Part (2) of the proposition
follows since $[I_{W \times_Y X} \to \cO_{W\times X}] = [\cO_{W
\times_Y X}]$ in the derived category.

We show the ``only if'' implication of part (3): suppose $Z\to W$ is
given such that the composite $Z \to Y$ satisfies the assumption in
part (3).  The pullback of $F = I_{W \times_Y X}$ to $Z\times X$ is
flat over $Z$, by the lemma it is torsion free, and as above it
injects to $\cO_{Z\times X}$, therefore it is an ideal sheaf,
which coincides with the ideal of $Z\times_Y X$ by right exactness of
the tensor product. 

Assume conversely that $I_{Z\times_Y X}$ is flat over $Z$. Its direct
image on $ Z \times Y$ is the ideal $I_{\Gamma_{Z \to Y}}$ of the graph
of $Z \to Y$, and by definition the inverse image of $I_{\Gamma_{Z \to
Y}}$ surjects onto $I_{Z\times_Y X}$. Thus $I_{Z\times_Y X}$ is a
perverse point-ideal sheaf, and therefore $Z \to Y$ factors through
$W$, proving the ``if'' implication of part (3).
\qed

\section{Surface contractions and Segre cones}

We now address Proposition \ref{Prop:surfaces}. 
\subsection{Smooth surface contractions} First let $f:X\to Y$ be
the contraction of smooth surfaces.  By the flatness criterion, to show
that $W  = Y$ it suffices to show that the ideal of the graph of $f$
in $X\times Y$ is flat over $Y$. The problem being local in the
\'etale topology, we may
assume $Y = \bbA^2$, with coordinates $x,y$, and $X = \bbA^2$,
with coordinates $u,v$, with $x=f_1(u,v)$ and $y = f_2(u,v)$. Then the ideal
$I$ is generated by 
$$ x-f_1(u,v), \ \ y-f_2(u,v).$$ The ideal is complete intersection, with
Koszul resolution 
$$\xymatrix{0\ar[r]&  \cO_{X\times Y} \ar[rrr]^{( y-f_2, f_1-x)}&&&
\cO_{X\times 
Y}^2\ar[rrr]^{\binom{x-f_1}{y-f_2}}&&& I \ar[r]&
0.}$$ 
We use the local criterion for flatness. To check that $I$ is flat over $Y$ it
suffices to check that 
$Tor_1^{\cO_Y}(I, E) = 0 $ where $E$ is the structure sheaf of
the center $P$ of blowing up, with maximal ideal $(x,y)$, and for this it
suffices to show that 

$$\xymatrix{ \cO_X \otimes E\ar[rr]^{( y-f_2, f_1-x)}&& (\cO_X \otimes E)^2}$$
is still injective. The 
latter coincides with the map 
$$\xymatrix{ \cO_X \ar[rr]^{( -f_2, f_1)}&& \cO_X^2},$$ which is
injective since 
$f_i$ is not identically 0. This
proves part (1) of Proposition \ref{Prop:surfaces}.

\subsection{Contractions of $(-2)$ curves} 
Now let $X \to Y$ be the contraction of a $(-2)$-curve. To show $W^n =
X$ it suffices to show that the 
\begin{enumerate}
\item[(a)] ideal of $X\times_YX$ is flat over
$X$, and 
\item[(b)] that the ideal of the graph of $f$ is {\em not} flat over $Y$.
\end{enumerate}

Again we can work \'etale locally, and reduce to the case where $Y =
\Spec \CC [x,y,z]/ (xy-z^2)$. one chart of the blow up is $X = \Spec \CC[u,v]$, with the
morphism given by 
$$x = u,\ \ z = u v, \ \ y = uv^2.$$  For the other  chart we have
$$y = u,\ \ z = u v, \ \ x = uv^2.$$  

To prove (a) we need to consider several affine charts on the
product. We describe one, leaving the others as an easy exercise to
the reader. Consider the affine chart on $X\times X$ with coordinates
$u_1,v_1,u_2,v_2$, where the ideal is 
$$I_{X\times_YX}\ =\ (\,u_1-u_2\,,\,u_1v_1-u_2v_2\,,\,u_1v_1^2-u_2v_2^2\,).$$
Note that the first two generators already generate the ideal, so 
 $I_{X\times_YX} = (u_1-u_2,\,u_1v_1-u_2v_2)$. The Koszul resolution
gives the result as in part (1). 

For (b), the ideal of the graph is generated by $$x-u,\ \ z-uv,\ \
y-uv^2.$$ It is a straightforward, though a bit tedious, to check that
it fails to be flat over $Y$. Alternatively, one can use the following
Macaulay 2 code:
\begin{verbatim}
        R=QQ[x,y,z,u,v]/(x*y-z^2);
        J=module(ideal(x-u, z-u*v,y-u*v^2));
        K = R^1/module(ideal(x,y,z));
        Tor_1(J,K)==0
\end{verbatim}

Of which the result is:
\begin{verbatim}
        o4 = false  
\end{verbatim}

This completes the proof of Proposition \ref{Prop:surfaces}.
\subsection{Other surface contractions}
 Proposition \ref{Prop:surfaces} leaves out many interesting cases of
 surface contractions, in particular the case of arbitrary crepant contractions
 of surfaces with only rational 
 double points.

The latter class of contractions turns out to follow from the work of the
second 
author \cite{Chen}. Indeed, any such crepant contraction locally
admits a one-parameter deformation deformation $\cX \to \cY\to C$,
so-called simultaneous partial 
resolution,   which on the total space is a threefold terminal
flopping contraction, where $\cX^+$ is the other simultaneous partial
resolution. In \cite{Chen} it is shown that $W(\cX / \cY) = 
\cX^+$. It is also shown that the central fiber of $W(\cX / \cY)$ is
$W(X/Y).$ Finally it is known that the fibers of the simultaneous partial
resolution $\cX^+$ are 
isomorphic to those of $\cX$, which shows that $W(X/Y) = X$.
\subsection{The Segre cone}
We now address Proposition \ref{Prop:segre}. Again, it suffices to
show that the ideal of $X\times_YX^+$ is flat over
$X^+$. Equivalently, it suffices to show that the structure sheaf of
$X\times_YX^+$ has $Tor$-dimension 1. We can work locally on $X^+$,
where an affine chart can be taken as $\bbA^{k+2}$, and the inverse
image of this chart in  $X\times_YX^+$ is the blowing up of $\bbA^k
\times \{0,0\}$. The fact that this has $Tor$-dimension 1 follows from
the case of a smooth surface (Proposition \ref{Prop:surfaces}, part
1), or can be seen directly: an affine chart for $X\times_YX^+$ is
 the hypersurface in $\bbA^{k+2}\times \bbA^1$, with coordinats
 $z_1,\ldots z_n, x,y,t$  given by equation $y=xt$. We thus have a
 resolution
$$0 \to \cO_{\bbA^{k+2}\times \bbA^1} \to  \cO_{\bbA^{k+2}\times
  \bbA^1} \to \cO_{X\times_YX^+} \to 0,$$
giving the result.

\section{The Francia flips}
We now address Proposition \ref{Prop:Francia} about the Francia flips
of index $n$. This is used  in the paper \cite{ABC}.
\subsection{The setup}\label{Sec:setup}
We follow the notation in Koll\'ar-Mori \cite{Kollar-Mori}, page 39:
denote by $Y_1$ the usual threefold quadratic 
singularity $xy = uv$, and let $X_1\to Y_1$ be one of the two small
resolutions, blowing up $(x,v)$, so $X_1 \to Y_1$ is the standard flopping
$(-1,-1)$ contraction. We denote by $X_1^+\to Y_1$ the other resolution, 
blowing up the ideal $(x,u)$. Denote
by $C_n$ the cyclic group of $n$ elements, acting on $Y_1$ via
$(x,y,u,v) \mapsto (\zeta x , y, \zeta u , v)$, which lifts to an
action on $X_1$ and on $X_1^+$.   The quotients of $X_1,Y_1,X_1^+$ are denoted
$X_n, Y_n, X^+_n$ respectively. 

In the proposition, we need to prove that $W^n(X_n/Y_n)$ is $X_n^+$.

\subsection{Flatness of the ideal}
By Corollary \ref{Cor:flp}  it suffices to verify that  the ideal of
$X_n \times_{Y_n}X^+_n$ is flat over $X_n^+$. 

The case $n=1$ follows from the main result of Bridgeland \cite{B}  (and can 
easily be checked by an explicit calculation). We thus have that the
ideal of $X_1 \times_{Y_1} X^+_1$ is flat over $X_1$. Since $X_n$ is
smooth and $X_1 \to X_n$ is finite, we have that $X_1 \to X_n$ is
flat, therefore  the
ideal of $X_1 \times_{Y_1} X^+_1$ is flat over $X_n$ as well.

We have an exact sequence of $\cO_{X_1 \times X^+_1}$-modules 
$$ 0 \to\ I_{X_1 \times_{Y_1} X^+_1}\ \to\  \cO_{X_1 \times X^+_1}\
\to \ \cO_{X_1 \times_{Y_1} X^+_1}\ \to 0. $$
Consider the diagonal action of $C_n$ on   $ X_1 \times X^+_1$. Taking
invariants in the exact sequence above we get a sequence
$$ 0 \to\ I_{(X_1 \times_{Y_1} X^+_1)/C_n}\ \to\  \cO_{(X_1 \times X^+_1)/C_n}\
\to \ \cO_{(X_1 \times_{Y_1} X^+_1)/C_n}\ \to 0, $$
which is still exact (since taking invariants of a finite group is exact in
characteristic 0).
Since the morphism $ X_1 \times X^+_1 \to X_1^+$ is equivariant, we
get an induced morphism $ (X_1 \times X^+_1)/C_n \to X_n^+$. Since
$I_{X_1 \times_{Y_1} X^+_1}$ is flat over $X_n^+$, and since the
module of  invariants $I_{(X_1 \times_{Y_1} X^+_1)/C_n}$ is a direct
summand, the latter is also flat over $X_n^+$. This means 
$$Tor_1(I_{(X_1 \times_{Y_1} X^+_1)/C_n}, M) = 0$$
fro every $\cO_{X_n^+}$-module $M$, or, equivalently,
$$Tor_2(\cO_{(X_1 \times_{Y_1} X^+_1)/C_n}, M) = 0.$$

The Lemma below says that $(X_1 \times_{Y_1} X^+_1)/C_n \to X_n
\times_{Y_n}X^+_n$ is an isomorphism. This means that 
$$Tor_2(\cO_{X_n \times_{Y_n}X^+_n}, M) = 0,$$ which in turn implies
that 
$$Tor_1(I_{X_n \times_{Y_n}X^+_n}, M) = 0,$$
so the ideal is flat \qed

We used the following:
\begin{lemma} The natural morphism $(X_1 \times_{Y_1} X^+_1)/C_n \to X_n
\times_{Y_n}X^+_n$ is an isomorphism
\end{lemma}

{\em Proof.} The morphism is finite and birational, so it suffices to
show that $X_n \times_{Y_n}X^+_n$ is normal. We check this by direct
computation.

We have $Y_1$ given by the equation $xy - uv$, therefore
 $$Y_n = \Spec \CC [\left\{x^n\right\},
 \left\{x^{n-1}u\right\},\ldots,\left\{xu^{n-1}\right\},
 \left\{u^n\right\}, y,v]/R$$ 
where $R$ is some ideal of relations. Here, as well as later, expressions in
braces are variables, and the expression inside the brace is the image in the
function field of $Y_1$

We have two affine charts for $X_1$:
$$X_1(1) = \Spec \CC \left[\left\{{\frac{x}{v}}\right\},y,v\right],\quad
X_1(2) = \Spec \CC \left[\left\{{\frac{v}{x}}\right\},x,u\right]$$
 giving the following charts on
$X_n$:

$$X_n(1) = \Spec \CC \left[\left\{\frac{x^n}{v^n}\right\},y,v\right],$$
$$X_n(2)\ \ = \ \ \Spec \CC \left[\left\{x^n\right\},
  \left\{x^{n-1}u\right\},\ldots,\left\{xu^{n-1}\right\}, 
  \left\{u^n\right\},y,v,\left\{{\frac{v^n}{x^n}}\right\}\right]\ /\
R'$$
where $$R'= 
\left(\ldots,\left\{{\frac{v^n}{x^n}}\right\} \left\{x^n\right\} - 
  v^n,\ldots,\left\{{\frac{v^n}{x^n}}\right\} \left\{u^n\right\} -
  y^n\right)$$  

We have two affine charts for $X_1^+$:
$$X_1^+(1) = \Spec \CC \left[\left\{{\frac{x}{u}}\right\},y,u\right],\quad 
X_1^+(2) = \Spec \CC \left[{\left\{\frac{u}{x}\right\}},x,v\right]$$
 giving the following charts on
$X_n^+$:
$$X_n^+(1) = \Spec \CC
\left[\left\{{\frac{x}{u}}\right\},y,\left\{u^n\right\}\right],\quad  
X_n^+(2) = \Spec \CC
\left[\left\{{\frac{u}{x}}\right\},\left\{x^n\right\},v\right].$$ 

We have the following charts on the fibered product:
\begin{eqnarray*}
X_n(1) \mathop\times\limits_{Y_n}X_n^+(1) & = & \Spec \CC
\left[\left\{{\frac{x^n}{v^n}}\right\},y,\left\{{\frac{x}{u}}\right\}\right] \\
X_n(1) \mathop\times\limits_{Y_n}X_n^+(2) & = & \Spec \CC
\left[\left\{{\frac{x^n}{v^n}}\right\},v,\left\{{\frac{u}{x}}\right\}\right] \\
X_n(2) \mathop\times\limits_{Y_n}X_n^+(1) & = & \Spec \CC
\left[\left\{{\frac{v^n}{x^n}}\right\},
\left\{u^n\right\},y,\left\{{\frac{x}{u}}\right\}\right]\ / \ 
\left(\left\{{\frac{v^n}{x^n}}\right\}   \left\{u^n\right\} - y^n\right) \\
X_n(2) \mathop\times\limits_{Y_n}X_n^+(2) & = & \Spec \CC
\left[\left\{{\frac{v^n}{x^n}}\right\},
\left\{x^n\right\},v,\left\{{\frac{u}{x}}\right\}\right]\ / \
\left(\left\{{\frac{v^n}{x^n}}\right\}   
  \left\{x^n\right\} - v^n\right) 
\end{eqnarray*}

The first two are smooth, the last two are a product of $\bbA^1$ with
an $A_{n-1}$-surface singularity, which is normal. \qed

\subsection{Non-flatness of the ideal}
The second part of the proposition boils down to showing that
$I_{X_2\times_{Y_2}X_2^+} $ is not flat over $X_2$. This can be done by a
direct calculation, e.g. the following Macaulay 2 code. The first output line
verifies the result already proven, namely flatness over $X_2^+$. The second
shows non-flatness over $X_2$:

\begin{verbatim}
R = QQ[E,G,H,A,B,C];
S = QQ[e,g,h];
F = map(S,R,{e^2,g^2,h^2,e*g,e*h,g*h});
T = QQ[a,b,c] ** image F ;
J = module(ideal(E-a^2*c,G-c,A-a*c,C-b,B-a*b));
K = T^1/module(ideal(a,b,c));
L = T^1/module(ideal(A,B,C,E,G,H));
Tor_1(J,K)==0
Tor_1(J,L)==0
\end{verbatim}
with output
\begin{verbatim}
o8 = true
o9 = false
\end{verbatim}

\end{document}